\documentclass[12pt]{article}
\usepackage{amsmath,amssymb,amsthm}
\usepackage[T1]{fontenc}
\usepackage[utf8]{inputenc}
\usepackage{epigraph}

\usepackage{float}
\usepackage[hyphens]{url}
\usepackage{authblk}       
\usepackage{hyperref}
\hypersetup{
  pdftitle  = {Beyond the Euler--Mascheroni Constant: A Family of Functionals},
  pdfauthor = {Ken Nagai}
}
\newcommand{\MSC}[1]{\par\small\textbf{MSC 2020.} #1\par}
\newcommand{\keywords}[1]{\par\small\textbf{Keywords.} #1\par}

\newtheorem{theorem}{Theorem}
\newtheorem{prop}{Proposition}
\newtheorem{corollary}{Corollary}
\newtheorem{remark}{Remark}

\newcommand{\Cl}{\operatorname{Cl}}

\title{Beyond the Euler--Mascheroni Constant:\\ A Family of Functionals}
\author{Ken Nagai\thanks{Email: \texttt{tknagai@outlook.com}. Independent Researcher.}}
\date{}

\begin{document}
\maketitle


\epigraph{\textit{``... in der Mathematik giebt es, kein Ignorabimus!''}}
{David Hilbert (1900)}

\begin{abstract}
We introduce a family of regularized functionals $g_n(x)$ that generalize
the Euler--Mascheroni constant $\gamma$.
They arise from a weighted regularization of Clausen-type trigonometric sums,
and admit explicit integral representations, differential and ladder relations,
together with an umbral generating function.
\end{abstract}

\MSC{11M06; 11B68; 33B15; 33E20}

\section{Introduction}

The Euler--Mascheroni constant $\gamma$ is classically defined as the finite part
of the divergent harmonic series
\cite{Apostol1976,GourdonSebah2008},
and it also appears as the constant term in the Laurent expansion of $\zeta(s)$
at $s=1$. However, in Clausen-type trigonometric sums
\cite{Lewin1981,CvijovicKlinowski1997},
the same divergence reappears in a more structured form.
We show in this note that by averaging with umbral weights $(1-t)^n$ and
introducing a scale $x$, the role of $\gamma$ is naturally generalized to an
infinite family of regularized functionals $g_n(x)$.
These functions admit integral representations, ladder relations,
and an umbral generating function
\cite{Roman1984,Blasiak2005,RotaTaylor1994}.

\begin{theorem}[Main Theorem]
The Euler--Mascheroni constant $\gamma$, arising as the finite part of
$\zeta(1)$, admits a natural extension to an infinite family of
regularized functionals $g_n(x)$.
They admit explicit integral representations, differential and
ladder relations, and a unified umbral generating function, thereby situating
$\gamma$ within a broader analytic framework.
\end{theorem}

\medskip
\noindent\textbf{Outline.}
Section~2 states our core identities, including an integral
representation, differential and ladder relations.
Section~3 presents the Bernoulli connections.
Section~4 develops the umbral generating function.
Section~5 concludes with a sketch of proof for the Main Theorem
and brief outlook remarks.

\section{Core identities}

\begin{prop}[Integral representation]
For $n\ge1$ and $x>0$,
\[
g_n(x)=H_n-\log(2\pi x)
- n \int_0^1 (1-u)^{n-1}\,\log\!\bigl(2\sin(\pi x u)\bigr)\,du.
\]
\end{prop}

\begin{corollary}[Small-$x$ normalization]
For $n\ge1$,
\[
g_n(x) = O(x^2)\qquad (x\to 0^+).
\]
\end{corollary}

\begin{prop}[Differential form]
For $n\ge1$ and $x>0$,
\begin{align*}
x\,g_n'(x)
&= 2n!\sum_{m\ge1}\frac{(2m)!}{(2m+n)!}\,\zeta(2m)\,x^{2m} - 1\\
&= -\,n \int_0^1 (1-u)^{n-1}\,\pi x u\,\cot(\pi x u)\,du \;-\; 1.
\end{align*}
\end{prop}

\begin{remark}[Cotangent connection]
The derivative formula shows that $g_n$ is governed by cotangent averages,
linking it with classical Bernoulli expansions of $\pi\cot(\pi z)$.
\end{remark}

\begin{prop}[Ladder in $n$]
For $n\ge1$ and $x>0$,
\[
g_{n+1}(x)-g_n(x)
=\frac{1}{n+1}
- \int_0^1\!\bigl[(n+1)(1-u)^n-n(1-u)^{n-1}\bigr]\,
\log(2\sin(\pi x u))\,du.
\]
\end{prop}

\begin{corollary}[Large-$n$ decay]
For fixed $x>0$,
\[
g_n(x) = O\!\left(\tfrac{1}{n}\right)\qquad (n\to\infty).
\]
\end{corollary}

\section{Bernoulli connections}

The cotangent expansion
\[
\pi \cot(\pi z) = \frac{1}{z} + \sum_{m=1}^\infty (-1)^m \frac{2^{2m} B_{2m}}{(2m)!}\,(\pi z)^{2m-1}
\]
implies that the integrands appearing in the differential form of $g_n(x)$
are closely tied to Bernoulli numbers $B_{2m}$.
Indeed, expanding $\pi x u \cot(\pi x u)$ and integrating termwise shows that
the coefficients in the defining series of $g_n(x)$ match precisely the
Bernoulli--zeta relation
\[
\zeta(2m) = (-1)^{m+1}\frac{(2\pi)^{2m}B_{2m}}{2(2m)!}.
\]

\begin{remark}[Bernoulli connection]
Thus $g_n(x)$ interpolates between the Euler--Mascheroni constant
and higher Bernoulli data, reflecting the classical bridge
between Clausen sums and even zeta values.
\end{remark}

\begin{remark}[Clausen origin]
In the cosine--Clausen expansion one encounters
\[
\Cl^{(c)}_{2m+1}(\theta)
= \sum_{r=0}^{m-1}\zeta(2m+1-2r)\,\frac{(-1)^r \theta^{2r}}{(2r)!}
 - \frac{(-1)^m \theta^{2m}}{(2m)!}\,g_{2m}\!\Bigl(\tfrac{|\theta|}{2\pi}\Bigr)
 + \cdots,
\]
where the divergent central contribution involving $\zeta(1)$
is naturally replaced by $g_{2m}(|\theta|/2\pi)$.
This indicates that the family $g_n(x)$ has its origin
in the same mechanism that produces the Euler--Mascheroni constant
in the simplest case $m=0$.
\end{remark}

\section{Umbral generating function}

\begin{prop}[Generating function]
For $|z|<1$,
\begin{align*}
\sum_{n\ge1} g_n(x)\,z^n
&= -\,\frac{z}{1-z}\,\log(2\pi x)
   - \frac{\log(1-z)}{1-z} \\
&\quad - \int_0^1 \log\!\bigl(2\sin(\pi x u)\bigr)\,
      \frac{z}{\bigl(1 - z(1-u)\bigr)^2}\,du.
\end{align*}
\end{prop}

\begin{remark}[Umbral perspective]
The generating function realizes $g_n(x)$
as coefficients of an operator-valued series
acting on $\log(2\sin(\pi x u))$.
This matches the umbral perspective
\cite{Roman1984,RotaTaylor1994,Blasiak2005},
where moments are encoded by a single parent functional.
\end{remark}

\section{Conclusion}

We close with a brief sketch of proof for the Main Theorem, 
together with some outlook remarks.

\begin{proof}[Sketch of proof]
The defining moment representation of $g_n(x)$
leads directly to the integral identity in Proposition~2.1.
Differentiating under the integral sign yields
the differential and ladder relations.
Finally, the umbral generating function provides a
unified operator perspective, from which the Main Theorem follows
by a straightforward synthesis.
\end{proof}

Beyond this, the umbral/Heisenberg perspective suggests
further connections with spectral theory and operator calculus.
In particular, the moment--functional interpretation points toward
an operator-theoretic framework in which $g_n(x)$
appears naturally alongside Bernoulli-type structures.
We leave such perspectives for future investigation.

\section*{Acknowledgments}
The author gratefully acknowledges the assistance of an AI language model (“fuga”) 
for valuable help with document structuring, stylistic polishing, and proofreading. Any remaining errors are my own.

\appendix
\section{Appendix: Numerical checks}

\begin{table}[H]
\centering
\begin{tabular}{c|c|c}
$(n,x)$ & Series definition & Integral representation \\
\hline
$(1,0.5)$ & $0.0770$ & $0.0770$ \\
$(2,0.5)$ & $-0.0619$ & $-0.0619$ \\
$(3,0.5)$ & $-0.0597$ & $-0.0597$ \\
$(2,1.0)$ & $-0.1639$ & $-0.1639$
\end{tabular}
\caption{Numerical consistency between the series and integral representations of $g_n(x)$.}
\end{table}

\begingroup
\raggedright\sloppy\setlength{\emergencystretch}{3em}

\end{document}